\documentclass[10pt,a4paper]{article}

\usepackage{amsmath}
\usepackage{amsthm}
\usepackage{amssymb}
\usepackage{amscd}

\theoremstyle{definition}
\newtheorem{theorem}{Theorem}[section]

\newtheorem{lemma}[theorem]{Lemma}

\newenvironment{demo}[1]{%
  \trivlist
  \item[\hskip\labelsep
        {\bf #1.}]
}{%
\hfill\qedsymbol
  \endtrivlist
}

\newcommand\Pf{\operatorname{Pf}}
\newcommand\perm{\operatorname{perm}}
\newcommand\Hf{\operatorname{Hf}}
\newcommand\Res{\operatorname{Res}}
\newcommand\Sym{{\mathcal{S}}}

\title{
A Pfaffian--Hafnian Analogue of Borchardt's Identity
}

\author{
Masao ISHIKAWA\thanks{Faculty of Education, Tottori University,
 e-mail:\tt{ishikawa@fed.tottori-u.ac.jp}}
\and
Hiroyuki KAWAMUKO\thanks{Faculty of Education, Mie University,
 e-mail:\tt{kawam@edu.mie-u.ac.jp}} 
\and
Soichi OKADA\thanks{Graduate School of Mathematics, Nagoya University,
 e-mail:\tt{okada@math.nagoya-u.ac.jp}}
}

\date{}

\begin{document}

\maketitle

\begin{abstract}
We prove
$$
\Pf \left(
 \frac{ x_i - x_j }{ (x_i + x_j)^2 }
 \right)_{1 \le i, j \le 2n}
 =
\prod_{1 \le i < j \le 2n} \frac{ x_i - x_j }{ x_i + x_j }
\cdot
\Hf \left(
 \frac{ 1 }{ x_i + x_j }
 \right)_{1 \le i, j \le 2n}
$$
(and its variants) by using the complex analysis.
This identity can be regarded as a Pfaffian--Hafnian analogue of
 Borchardt's identity and as a generalization of
 Schur's identity.
\end{abstract}

\section{
Introduction
}

Determinant and Pfaffian identities play a key role in combinatorics
 and the representation theory (see, for example, \cite{IOW}, \cite{IW1},
 \cite{IW2}, \cite{Kr}, \cite{O1}, \cite{O2}).
Among such determinant identities, the central ones are
 Cauchy's determinant identities
 (\cite{C})
\begin{gather}
\det \left(
 \frac{ 1 }{ x_i + y_j }
\right)_{1 \le i, j \le n}
 =
\frac{ \prod_{1 \le i < j \le n} (x_j - x_i) (y_j - y_i) }
     { \prod_{i,j=1}^n (x_i + y_j) },
\label{cauchy1}
\\
\det \left(
 \frac{ 1 }{ 1 - x_i y_j }
\right)_{1 \le i, j \le n}
 =
\frac{ \prod_{1 \le i < j \le n} (x_j - x_i) (y_j - y_i) }
     { \prod_{i,j=1}^n (1 - x_i y_j) }.
\label{cauchy2}
\end{gather}
C.~W.~Borchardt \cite{B} gave a generalization of Cauchy's identities:
\begin{gather}
\det \left(
 \frac{ 1 }{ (x_i + y_j)^2 }
\right)_{1 \le i, j \le n}
 =
\frac{ \prod_{1 \le i < j \le n} (x_j - x_i) (y_j - y_i) }
     { \prod_{i,j=1}^n (x_i + y_j) }
\cdot
\perm \left(
 \frac{ 1 }{ x_i + y_j }
\right)_{1 \le i, j \le n},
\label{borch1}
\\
\det \left(
 \frac{ 1 }{ (1 - x_i y_j)^2 }
\right)_{1 \le i, j \le n}
 =
\frac{ \prod_{1 \le i < j \le n} (x_j - x_i) (y_j - y_i) }
     { \prod_{i,j=1}^n (1 - x_i y_j) }
\cdot
\perm \left(
 \frac{ 1 }{ 1 - x_i y_j }
\right)_{1 \le i, j \le n}.
\label{borch2}
\end{gather}
Here $\perm A$ is the permanent of a square matrix $A$ defined by
$$
\perm A =
 \sum_{\sigma \in \Sym_n}
  a_{1 \sigma(1)} a_{2 \sigma(2)} \cdots a_{n \sigma(n)}.
$$
This identity (\ref{borch1}) is used when we evaluate the determinants
 appearing in the $0$-enumeration of alternating sign matrices
 (see \cite{O2}).

I.~Schur \cite{S} gave a Pfaffian analogue of Cauchy's identity
 (\ref{cauchy1})
 in his study of projective representations of the symmetric groups.
Schur's Pfaffian identity and its variant (\cite{LLT}, \cite{St})
 are
\begin{gather}
\Pf \left( \frac{x_j - x_i}{x_j + x_i} \right)_{1 \le i, j \le 2n}
 =
\prod_{1 \le i < j \le 2n} \frac{x_j - x_i}{x_j+x_i},
\label{schur1}
\\
\Pf \left( \frac{x_j - x_i}{1 - x_i x_j} \right)_{1 \le i, j \le 2n}
 =
\prod_{1 \le i < j \le 2n} \frac{x_j - x_i}{1 - x_i x_j}.
\label{schur2}
\end{gather}

In this note, we give identities which can be regarded
 as Pfaffian analogues of Borchardt's identities
 (\ref{borch1}), (\ref{borch2})
 and as generalizations of Schur's identities
 (\ref{schur1}), (\ref{schur2}).

\begin{theorem}
\label{main-theorem}
Let $n$ be a positive integer.
Then we have
\begin{gather}
\Pf \left(
 \frac{ x_i - x_j }{ (x_i + x_j)^2 }
 \right)_{1 \le i, j \le 2n}
 =
\prod_{1 \le i < j \le 2n} \frac{ x_i - x_j }{ x_i + x_j }
\cdot
\Hf \left(
 \frac{ 1 }{ x_i + x_j }
 \right)_{1 \le i, j \le 2n},
\label{main1}
\\
\Pf \left(
 \frac{ x_i - x_j }{ (1 - x_i x_j)^2 }
 \right)_{1 \le i, j \le 2n}
 =
\prod_{1 \le i < j \le 2n} \frac{ x_i - x_j }{ 1 - x_i x_j }
\cdot
\Hf \left(
 \frac{ 1 }{ 1 - x_i x_j }
 \right)_{1 \le i, j \le 2n}.
\label{main2}
\end{gather}
\end{theorem}

Here $\Hf A$ denotes the Hafnian of a symmetric matrix $A$ defined by
$$
\Hf A = 
 \sum_{\sigma \in {\mathcal{F}}_{2n}}
  a_{\sigma(1)\sigma(2)} a_{\sigma(3)\sigma(4)} \cdots
  a_{\sigma(2n-1)\sigma(2n)},
$$
where $\mathcal{F}_{2n}$ is the set of all permutations $\sigma$ satisfying
 $\sigma(1) < \sigma(3) < \cdots < \sigma(2n-1)$ and
 $\sigma(2i-1) < \sigma(2i)$ for $1 \le i \le n$.

\section{
Proof
}

In this section, we prove the identity (\ref{main1}) in
 Theorem~\ref{main-theorem} by using the complex analysis.
The other identity (\ref{main2}) is shown by the same method,
 and also derived from more general identity (\ref{main3}) in
 Theorem~\ref{general2}, which follows from (\ref{main1}).
So we omit the proof of (\ref{main2}) here.

Hereafter we put
$$
A = \left( \frac{x_i - x_j}{(x_i+x_j)^2} \right)_{1 \le i, j \le 2n},
\quad
B = \left( \frac{1}{x_i + x_j} \right)_{1 \le i, j \le 2n}.
$$
For an $2n \times 2n$ symmetric (or skew-symmetric) matrix $M = (m_{ij})$
 and distinct indices $i_1, \cdots, i_r$, we denote by
 $M^{i_1, \cdots, i_r}$ the $(2n-r) \times (2n-r)$ matrix obtained by
 removing the rows and columns indexed by $i_1, \cdots, i_r$.

First we show two lemmas by using the complex analysis.

\begin{lemma}
\label{lemma1}
\begin{equation}
\sum_{\substack{1 \le k, l \le 2n \\ k \neq l}}
 \frac{1}{(x_k - z)(x_l + z)} \Hf(B^{k,l})
 =
\Hf(B) \cdot \sum_{k=1}^{2n} \frac{2x_k}{x_k^2 - z^2}.
\label{eq-lemma1}
\end{equation}
\end{lemma}

\begin{demo}{Proof}
Let us denote by $F(z)$ (resp. $G(z)$) the left (resp. right) hand side
 of (\ref{eq-lemma1}),
 and regard $F(z)$ and $G(z)$ as rational functions
 in the complex variable $z$,
 where $x_1, \cdots, x_{2n}$ are distinct complex numbers.
Then $F(z)$ and $G(z)$ have poles at $z = \pm x_1, \cdots, \pm x_{2n}$
 of order $1$.
The residues of $F(z)$ at $z = \pm x_m$ are given by
$$
\Res_{z = x_m} F(z)
 =
- \sum_{\substack{1 \le l \le 2n \\ l \neq m}}
  \frac{1}{x_l + x_m} \Hf (B^{m,l}),
\quad
\Res_{z = -x_m} F(z)
 =
\sum_{\substack{1 \le k \le 2n \\ k \neq m}}
  \frac{1}{x_k + x_m} \Hf( B^{k,m}).
$$
By considering the expansion of $\Hf (B)$ along the $m$th row/column,
 we have
$$
\Res_{z = x_m} F(z) = - \Hf (B),
\quad
\Res_{z = -x_m} F(z) = \Hf (B).
$$
On the other hand, the residues of $G(z)$ at $z = \pm x_m$ are given by
\begin{gather*}
\Res_{z = x_m} G(z)
 = -\Hf(B) \cdot \frac{2x_m}{2x_m}
 = -\Hf(B),
\\
\Res_{z = -x_m} G(z)
 = \Hf(B) \cdot \frac{2x_m}{2x_m}
 = \Hf (B).
\end{gather*}
Since $\lim_{z \to \infty} F(z) = \lim_{z \to \infty} G(z) = 0$,
 we conclude that $F(z) = G(z).$
\end{demo}

\begin{lemma}
\label{lemma2}
If $n$ is a positive integer,
then
\begin{equation}
\sum_{k=1}^{2n-1}
 \frac{x_k - z}{(x_k + z)^2}
 \prod_{\substack{1 \le i \le 2n-1 \\ i \neq k}} \frac{ x_k + x_i }{x_k - x_i}
 \cdot
 \Hf (B^{k,2n})
 =
 \prod_{i=1}^{2n-1} \frac{x_i - z}{x_i + z}
 \sum_{k=1}^{2n-1} \frac{1}{x_k + z} \Hf(B^{k,2n}).
\label{eq-lemma2}
\end{equation}
\end{lemma}

\begin{demo}{Proof}
Let $P(z)$ (resp. $Q(z)$) be the left (resp. right) hand side
 of (\ref{eq-lemma2}),
 and regard $P(z)$ and $Q(z)$ as rational functions in $z$,
 where $x_1, \cdots, x_{2n-1}$ are distinct complex numbers.
Then $P(z)$ and $Q(z)$ have poles at $z = -x_1, \cdots, -x_{2n-1}$
 of order $2$.
Thus,
for a fixed $m$ such that $1\leq m\leq2n-1$,
we can write
\begin{align*}
P(z)
 &= \frac{p_2}{(z + x_m)^2} + \frac{p_1}{z + x_m}
 + O(z + x_m),
\\
Q(z)
 &= \frac{q_2}{(z + x_m )^2} + \frac{q_1}{z + x_m }
 + O(z + x_m),
\end{align*}
in a neighborhood of $z=-x_m$.
Now we compute the coefficients $p_2$, $p_1$, $q_2$ and $q_1$,
 and prove $p_2 = q_2$, $p_1 = q_1$.

By using the relation
$$
\frac{x_m - z}{(x_m + z)^2}
 = \frac{2x_m}{(x_m + z)^2} - \frac{1}{x_m + z},
$$
we see that
\begin{align}
p_2
 &=
 2x_m
 \prod_{\substack{1 \le i \le 2n-1 \\ i \neq m}} \frac{x_m + x_i}{x_m - x_i}
 \cdot \Hf (B^{m,2n}),
\label{p2}
\\
p_1
 &=
 - \prod_{\substack{1 \le i \le 2n-1 \\ i \neq m}} \frac{x_m + x_i}{x_m - x_i}
 \cdot \Hf (B^{m,2n}).
\label{p1}
\end{align}

Next we deal with
$$
Q(z)
 =
\frac{x_m - z}{x_m + z}
 \times
 \prod_{\substack{1 \le i \le 2n-1 \\ i \neq m}} \frac{x_i - z}{x_i + z}
 \times
 \sum_{k=1}^{2n-1} \frac{1}{x_k + z} \Hf (B^{k,2n}).
$$
The first factor can be written in the form
$$
\frac{x_m - z}{x_m + z}
 = \frac{2x_m}{x_m + z} - 1.
$$
By using the Taylor expansion $\log(1-t) = - t + O(t^2)$, we have
\begin{align*}
\log \frac{x_i - z}{x_m + x_i}
 &= - \frac{z + x_m}{x_i + x_m} + O \left( (z + x_m)^2 \right),
\\
\log \frac{x_i + z}{x_m - x_i}
 &= \frac{z + x_m}{x_i - x_m} + O \left( (z + x_m)^2 \right).
\end{align*}
Hence we see that
$$
\log \left(
 \frac{x_i - z}{x_i + z} \Bigm/ \frac{x_i + x_m}{x_i - x_m}
 \right)
=
 - \frac{2x_i}{x_i^2 - x_m^2}(z + x_m)
 + O\left( (z + x_m)^2 \right).
$$
Therefore the second factor of $Q(z)$ has the form
$$
\prod_{\substack{1 \le i \le 2n-1 \\ i \neq m}} \frac{x_i - z}{x_i + z}
 =
 \prod_{\substack{1 \le i \le 2n-1 \\ i \neq m}} \frac{x_i + x_m}{x_i - x_m}
 \cdot
 \left\{
  1
  - \sum_{\substack{1 \le k \le 2n-1 \\ k \neq m}} \frac{2x_k}{x_k^2 -x_m^2}
  \cdot (z + x_m)
  + O \left( (z + x_m)^2 \right)
 \right\}.
$$
Since we have
$$
\frac{1}{x_k + z} = \frac{1}{x_k - x_m} + O \left( z + x_m \right),
$$
the last factor of $Q(z)$ has the following expansion:
$$
\sum_{k=1}^{2n-1} \frac{1}{x_k + z} \Hf( B^{k,2n})
 = \frac{1}{x_m + z} \Hf (B^{m,2n})
 + \sum_{\substack{1 \le k \le 2n-1 \\ k \neq m}}
     \frac{1}{x_k - x_m} \Hf (B^{k,2n})
 + O( z + x_m ).
$$
Combining these expansions, we have
\begin{equation}
q_2
 =
 2 x_m
 \prod_{\substack{1 \le i \le 2n-1 \\ i \neq m}} \frac{x_i + x_m}{x_i - x_m}
 \cdot \Hf (B^{m,2n}),
\label{q2}
\end{equation}
and
\begin{align}
q_1
 &=
 \prod_{\substack{1 \le i \le 2n-1 \\ i \neq m}} \frac{x_i + x_m}{x_i - x_m}
\notag
\\
 &\quad\times
 \left\{
 2x_m \sum_{\substack{1 \le k \le 2n-1 \\ k \neq m}}
  \frac{\Hf (B^{k,2n})}{x_k - x_m}
  -
 2 x_m \Hf(B^{m,2n})
 \sum_{\substack{1 \le k \le 2n-1 \\ k \neq m}}
  \frac{2x_k}{x_k^2 - x_m^2}
 -
 \Hf (B^{m,2n})
 \right\}.
\label{q1}
\end{align}

It follows from (\ref{p2}) and (\ref{q2}) that $p_2 = q_2$.
From (\ref{p1}) and (\ref{q1}),
 in order to prove the equality $p_1 = q_1$,
it is enough to show that
$$
\sum_{\substack{1 \le k \le 2n-1 \\ k \neq m}}
 \frac{1}{x_k - x_m} \Hf(B^{k,2n})
 =
 \Hf (B^{m,2n})
 \sum_{\substack{1 \le k \le 2n-1 \\ k \neq m}} \frac{2x_k}{x_k^2 - x_m^2}.
$$
By permuting the variables $x_1, \cdots, x_{2n-1}$, we may assume that $m = 2n-1$.
Then, by expanding the Hafnian on the left hand side along the last
 row/column, it is enough to show that
$$
\sum_{k=1}^{2n-2} \frac{1}{x_k - x_{2n-1}}
 \sum_{\substack{1 \le l \le 2n-2 \\ l \neq k}}
  \frac{1}{x_l + x_{2n-1}} \Hf( B^{k,l,2n-1,n})
 =
\Hf (B^{2n-1,2n})
 \sum_{k=1}^{2n-2} \frac{2x_k}{x_k^2 - x_{2n-1}^2}.
$$
This follows from Lemma~\ref{lemma1} (with $2n$ replaced by $2n-2$ and
 $z$ replaced by $x_{2n-1}$),
 and we complete the proof of Lemma~\ref{lemma2}.
\end{demo}

Now we are in the position to prove the identity (\ref{main1})
 in Theorem~\ref{main-theorem}.

\begin{demo}{Proof of (\ref{main1})}
We proceed by induction on $n$.

Expanding the Pfaffian along the last row/column and
 using the induction hypothesis, we see
\begin{align*}
\Pf(A)
 &
 = \sum_{k=1}^{2n-1} (-1)^{k-1} \frac{x_k - x_{2n}}{(x_k + x_{2n})^2} \Pf (A^{k,2n})
\\
 &
 = \sum_{k=1}^{2n-1} (-1)^{k-1} \frac{x_k - x_{2n}}{(x_k + x_{2n})^2}
 \prod_{\substack{1 \le i < j \le 2n-1 \\ i,j \neq k}}
  \frac{x_i- x_j}{x_i + x_j}
  \Hf (B^{k,2n}).
\end{align*}
By using the relation
$$
\prod_{\substack{1 \le i < j \le 2n-1 \\ i,j \neq k}}
 \frac{x_i- x_j}{x_i + x_j}
 =
 (-1)^{k-1} \prod_{1 \le i < j \le 2n-1} \frac{x_i - x_j}{x_i + x_j}
 \cdot
 \prod_{\substack{1 \le i \le 2n-1 \\ i \neq k}} \frac{ x_k + x_i}{x_k - x_i},
$$
we have
$$
\Pf (A)
 =
\prod_{1 \le i < j \le 2n-1} \frac{x_i - x_j}{x_i + x_j}
\sum_{k=1}^{2n-1}
 \frac{x_k - x_{2n}}{(x_k + x_{2n})^2}
 \prod_{\substack{1 \le i \le 2n-1 \\ i \neq k}} \frac{x_k + x_i}{x_k - x_i}
 \cdot \Hf (B^{k,2n}).
$$
On the other hand, by expanding the Hafnian along the last row/column,
 we have
$$
\prod_{1 \le i < j \le 2n} \frac{ x_i - x_j }{  x_i + x_j }\cdot \Hf (B)
 =
\prod_{1 \le i < j \le 2n} \frac{x_i - x_j}{x_i + x_j}
 \sum_{k=1}^{2n-1} \frac{1}{x_k + x_{2n}} \cdot \Hf (B^{k,2n}).
$$
So it is enough to show the following identity:
$$
\sum_{k=1}^{2n-1}
 \frac{x_k - x_{2n}}{( x_k + x_{2n})^2}
 \prod_{\substack{1 \le i \le 2n-1 \\ i \neq k}} \frac{x_k + x_i}{x_k - x_i}
 \cdot\Hf (B^{k,2n})
 =
 \prod_{i=1}^{2n-1} \frac{x_i - x_{2n}}{x_i + x_{2n}}
 \sum_{k=1}^{2n-1} \frac{1}{x_k + x_{2n}} \cdot\Hf (B^{k,2n}).
$$
This identity follows from Lemma~\ref{lemma2} and the proof completes.
\end{demo}

\section{
Generalization
}

The Cauchy's identities (\ref{cauchy1}) and (\ref{cauchy2}),
 and the Borchardt's identities (\ref{borch1}) and (\ref{borch2})
 are respectively unified in the following form.

\begin{theorem}
\label{general1}
Let $f(x,y) = a x y + b x + c y + d$ be a nonzero polynomial.
Then we have
\begin{align}
&
\det \left( \frac{1}{f(x_i,y_j)} \right)_{1 \le i, j \le n}
\notag
\\
&\quad=
(-1)^{n(n-1)} (ad-bc)^{n(n-1)/2}
\frac{\prod_{1 \le i < j \le n} (x_j - x_i) (y_j - y_i) }
     {\prod_{1 \le i, j \le n} f(x_i,y_j) },
\label{cauchy3}
\\
&
\det \left( \frac{1}{f(x_i,y_j)^2} \right)_{1 \le i, j \le n}
\notag
\\
&\quad=
(-1)^{n(n-1)} (ad-bc)^{n(n-1)/2}
\frac{\prod_{1 \le i < j \le n} (x_j - x_i) (y_j - y_i) }
     {\prod_{1 \le i, j \le n} f(x_i,y_j) }
 \cdot
 \perm \left( \frac{1}{f(x_i,y_j)} \right)_{1 \le i, j \le n}.
\label{borch3}
\end{align}
\end{theorem}

Similarly we can generalize the Schur's identities (\ref{schur1}) and
 (\ref{schur2}),
 and our identities (\ref{main1}) and (\ref{main2}).

\begin{theorem}
\label{general2}
Let $g(x,y) = a x y + b (x + y) + c$ be a nonzero polynomial.
Then we have
\begin{align}
\Pf \left(
 \frac{ x_j - x_i }
      { g(x_i,x_j) }
\right)_{1 \le i, j \le 2n}
&=
(b^2 - ac)^{n(n-1)}
\prod_{1 \le i < j \le 2n}
 \frac{ x_j - x_i }
      { g(x_i,x_j) },
\label{schur3}
\\
\Pf \left(
 \frac{ x_j - x_i }
      { g(x_i,x_j)^2 }
\right)_{1 \le i, j \le 2n}
&=
(b^2 - ac)^{n(n-1)}
\prod_{1 \le i < j \le 2n}
 \frac{ x_j - x_i }
      { g(x_i,x_j) }
\Hf \left(
 \frac{ 1 }
      { g(x_i,x_j) }
\right)_{1 \le i, j \le 2n}.
\label{main3}
\end{align}
\end{theorem}

This generalization (\ref{schur3}) is given in \cite{Kn}.

\begin{demo}{Proof}
We derive (\ref{schur3}) and (\ref{main3}) from
 (\ref{schur1}) and (\ref{main1}) respectively.

First we consider the case where $b^2 - ac \neq 0$.
Suppose that $a \neq 0$.
Then, by putting
$$
A = \frac{1}{2},
\quad
B = \frac{1}{2a} ( b + \sqrt{b^2 - ac} ),
\quad
C = a,
\quad
D = b - \sqrt{b^2 - ac},
$$
and substituting
$$
x_i \to \frac{A x_i + B}{C x_i + D}
\quad(1 \le i \le 2n)
$$
in (\ref{schur1}) and (\ref{main1}),
 we obtain (\ref{schur3}) and (\ref{main3}).
Similarly we can show the case where $c \neq 0$.

If $b^2 - ac = 0 $ and $a \neq 0$, then we have
$$
g(x_i,x_j) = a^{-1} (a x_i + b)(a x_j + b).
$$
Hence we can evaluate the left hand sides of (\ref{schur3}) and (\ref{main3})
 by using
$$
\Pf \left( x_j - x_i \right)_{1 \le i, j \le 2n}
 =
\begin{cases}
 x_2 - x_1 &\text{if $n = 1$,} \\
 0         &\text{if $n \ge 2$,}
\end{cases}
$$
and obtain the equalities in (\ref{schur3}) and (\ref{main3}).
Similarly we can show the case where $b^2 - ac = 0$ and $c \neq 0$.
\end{demo}

From (\ref{cauchy3}) and (\ref{borch3}), we have
$$
\det \left( \frac{1}{f(x_i,y_j)^2} \right)_{1 \le i, j \le n}
 =
\det \left( \frac{1}{f(x_i,y_j)} \right)_{1 \le i, j \le n}
 \cdot
\perm \left( \frac{1}{f(x_i,y_j)} \right)_{1 \le i, j \le n}.
$$
Since the matrix $\left( f(x_i,y_j) \right)_{1 \le i, j \le n}$ has
 rank at most $2$, this identity is the special case of the following
 theorem.

\begin{theorem} \label{carlitz}
 (Carlitz and Levine \cite{CL})
Let $A = (a_{ij})$ be a matrix of rank at most $2$.
If $a_{ij} \neq 0$ for all $i$ and $j$, we have
$$
\det \left( \frac{1}{a_{ij}^2} \right)_{1 \le i, j \le n}
 =
\det \left( \frac{1}{a_{ij}} \right)_{1 \le i, j \le n}
 \cdot
\perm \left( \frac{1}{a_{ij}} \right)_{1 \le i, j \le n}.
$$
\end{theorem}

From (\ref{schur3}) and (\ref{main3}), we have
$$
\Pf \left(
 \frac{ x_j - x_i }
      { g(x_i,x_j)^2 }
\right)_{1 \le i, j \le 2n}
 =
\Pf \left(
 \frac{ x_j - x_i }
      { g(x_i,x_j) }
\right)_{1 \le i, j \le 2n}
 \cdot
\Hf \left(
 \frac{ 1 }
      { g(x_i,x_j) }
\right)_{1 \le i, j \le 2n}.
$$
It is a natural problem to find a Pfaffian--Hafnian analogue
 of Theorem~\ref{carlitz}.
Also it is interesting to find more examples of a skew-symmetric matrix $X$
 and a symmetric matrix $Y$ satisfying
$$
\Pf \left( x_{ij} y_{ij} \right)_{1 \le i, j \le 2n}
 =
\Pf \left( x_{ij} \right)_{1 \le i, j \le 2n}
 \cdot
\Hf \left( y_{ij} \right)_{1 \le i, j \le 2n}.
$$
Recently there appeared a bijective proof of Borchardt's identity
(see \cite{Si}).
It will be an interesting problem to give a bijective proof of \thetag{\ref{main1}} and \thetag{\ref{main2}}.


\begin{thebibliography}{99}

\bibitem{B}
C.~W.~Borchardt,
Bestimmung der symmetrischen Verbindungen vermittelst ihrer erzeugenden
 Funktion,
J. Reine Angew. Math. {\bf 53} (1855), 193--198.

\bibitem{C}
A.~L.~Cauchy,
M\'emoire sur les fonctions altert\'ees et sur les sommes altern\'ees,
Exercices Anal. et Phys. Math. {\bf 2} (1841), 151--159.

\bibitem{CL}
L.~Carlitz and J.~Levine,
An identity of Cayley,
Amer. Math. Monthly {\bf 67} (1960), 571--573.




\bibitem{IOW}
M.~Ishikawa, S.~Okada and M.~Wakayama,
Applications of minor summation formulas I : Littlewood's formulas,
J. Algebra {\bf 183} (1996), 193--216.

\bibitem{IW1}
M.~Ishikawa and M.~Wakayama,
Applications of minor summation formulas II : Pfaffians and Schur polynomials,
J. Combin. Theory Ser. A {\bf 88} (1999), 136--157.

\bibitem{IW2}
M.~Ishikawa and M.~Wakayama,
Applications of minor summation formulas III : Pl"ucker relations,
 lattice paths and Pfaffian identities,
{\tt arXiv:math.CO/0312358}.

\bibitem{Kn}
D.~Knuth,
Overlapping Pfaffians,
Electron. J. Combin. {\bf 3} (2) (``The Foata Festschrift'') (1996), 151--163.

\bibitem{Kr}
C.~Krattenthaler,
Advanced determinant calculus,
Sem. Lothar. Combin.
{\bf 42} (``The Andrews Festschrift'') (1999), Article B42q.

\bibitem{LLT}
D.~Laksov, A.~Lascoux and A.~Thorup,
On Giambelli's theorem on complete correlations,
Acta Math. {\bf 162} (1989), 143--199.

\bibitem{O1}
S.~Okada,
Application of minor summation formulas to rectangular-shaped
 representations of classical groups,
J. Algebra {\bf 205} (1998), 337--367.

\bibitem{O2}
S.~Okada,
Enumeration of symmetry classes of alternating sign matrices
 and characters of classical groups,
{\tt arXiv:math.CO/0308234},
to appear.

\bibitem{S}
I.~Schur,
\"Uber die Darstellung der symmetrischen und der alternirenden Gruppe
 durch gebrochene lineare Substitutuionen,
J. Reine Angew. Math. {\bf 139} (1911), 155--250.

\bibitem{Si}
D.~Singer,
A bijective proof of Borchardt's identity,
Electron. J. Combin. {\bf 11} (1), R48

\bibitem{St}
J.~R.~Stembridge,
Non-intersecting paths, Pfaffians and plane partitions,
Adv. Math. {\bf 83} (1990), 96--131.

\end{thebibliography}
\end{document}